\documentclass[12pt,leqno]{amsart}

\usepackage{amssymb}
 \usepackage{amsthm}

\usepackage{amsmath,amsfonts,amsbsy,latexsym,multicol,fancybox}
\usepackage{graphicx}
\usepackage{epstopdf}
\usepackage{color}
\DeclareGraphicsExtensions{.pdf, .jpg, .pnp, .bmp, .fig}
\usepackage{hyperref}

\usepackage{caption}\usepackage{subcaption}

\newtheorem{theorem}{Theorem}
\newtheorem{assumption}{Assumption}
\newtheorem{definition}{Definition}
\newtheorem{lemma}{Lemma}

\newcommand{\bass}{\begin{assumption}}\newcommand{\eass}{\end{assumption}}
\newcommand{\bde}{\begin{definition}} \newcommand{\ede}{\end{definition}}
\newcommand{\ble}{\begin{lemma}} \newcommand{\ele}{\end{lemma}}
\newcommand{\bth}{\begin{theorem}} \newcommand{\ethe}{\end{theorem}}

\newcommand{\bpf}{\begin{proof}}\newcommand{\epf}{\end{proof}}
\newcommand{\barr}{\begin{array}}\newcommand{\earr}{\end{array}}
\newcommand{\beao}{\begin{eqnarray*}}\newcommand{\eeao}{\end{eqnarray*}\noindent}
\newcommand{\beam}{\begin{eqnarray}}\newcommand{\eeam}{\end{eqnarray}\noindent}
\newcommand{\beqq}{\begin{equation}}\newcommand{\eeqq}{\end{equation}\noindent}

\newcommand{\ov}{\overline} 
\newcommand{\wt}{\widetilde}

\newcommand{\nto}{n\to\infty}

\newcommand{\tto}{t\to\infty}

\newcommand{\Rto}{R\to\infty}

\newcommand{\ga}{\gamma} 

\newcommand{\D}{\Delta}
  \newcommand{\ep}{\epsilon}
\newcommand{\ka}{\kappa}

\newcommand{\w}{\omega} \newcommand{\W}{\Omega}

\newcommand{\bfE}{{\mathbb E}}\newcommand{\bbE}{{\mathcal E}} 
\newcommand{\bbf}{{\mathcal F}}

\newcommand{\bbi}{{\mathbb I}}

\newcommand{\bbl}{{\mathcal L}}

 \newcommand{\bbN}{{\mathbb N}}

\newcommand{\bfP}{{\mathbb P}}
                              
 \newcommand{\bbR}{{\mathbb R}}

\graphicspath{{./Images/}}

\begin{document}




\title[A note on the asymptotic stability of the Semi-Discrete method for SDEs]{A note on the asymptotic stability of the Semi-Discrete method for Stochastic Differential Equations}
 \author{I. S. Stamatiou}
  \address{University of West Attica, Department of Biomedical Sciences}
 \email{joniou@gmail.com, istamatiou@uniwa.gr}
 \author{N. Halidias}
 \address{University of the Aegean, Department of Statistics and Actuarial-Financial Mathematics}
 \email{nick@aegean.gr}

\begin{abstract}
We study the asymptotic stability of the semi-discrete (SD) numerical method for the approximation of stochastic differential equations. Recently, we examined the order of $\bbl^2$-convergence of the truncated SD method and showed that it can be arbitrarily close to $1/2,$ see \textit{Stamatiou, Halidias (2019), Convergence rates of the Semi-Discrete method for stochastic differential equations, Theory of Stochastic Processes, 24(40)}. We show that the truncated SD method is able to preserve the asymptotic stability of the underlying SDE. Motivated by a numerical example, we also propose a different SD scheme, using the Lamperti transformation to the original SDE, which we call Lamperti semi-discrete (LSD). Numerical simulations support our theoretical findings. 
\end{abstract}

\maketitle
\tableofcontents

\section{Introduction}\label{asTSD:sec:intro}
\setcounter{equation}{0}

We study the following class of scalar stochastic differential equations (SDEs),
\beqq  \label{asTSD-eq:scalarSDEs}
dx_t =a(t,x_t)dt + b(t,x_t)dW_t, \qquad t\in[0,T],
\eeqq
where $a, b: [0,T]\times\bbR\rightarrow\bbR$ are measurable functions such that (\ref{asTSD-eq:scalarSDEs}) has a unique solution and $x_0$ is independent of all $\{W_{t}\}_{t\geq0}.$ We assume that SDE (\ref{asTSD-eq:scalarSDEs}) has non-autonomous coefficients, i.e. $a(t,x), b(t,x)$ depend explicitly on $t.$
SDEs of the type (\ref{asTSD-eq:scalarSDEs}) rarely have explicit solutions, therefore the need for numerical approximations for simulations of the solution process $x_t(\w)$ is apparent. In the case of nonlinear drift and diffusion coefficients classical methods may fail to strongly approximate (in the mean-square sense) the solution of (\ref{asTSD-eq:scalarSDEs}), c.f. \cite{hutzenthaler_et_al.:2011}, where the Euler method may explode in finite time. 

In this direction, we study the semi-discrete (SD) method originally proposed in \cite{halidias:2012} and further investigated  in  \cite{halidias_stamatiou:2016}, \cite{halidias:2014}, \cite{halidias:2015}, \cite{halidias:2015d}, \cite{halidias_stamatiou:2015} and recently in \cite{stamatiou:2018} and \cite{STAMATIOU:2019}. The main idea behind the semi-discrete method is freezing on each subinterval appropriate parts of the drift and diffusion coefficients of the solution at the beginning of the subinterval so as to obtain explicitly solved SDEs. Of course the way of freezing (discretization) is not unique.

The SD method is a fixed-time step explicit numerical method which strongly converges to the exact solution and also preserves the domain of the solution; if for instance the solution process $x_t$ is nonnegative then the approximation process $y_t$ is also nonnegative. The $\bbl^2$-convergence of the truncated SD method, see \cite{stamatiou_halidias:2019}, was recently shown to be arbitrarily close to $1/2.$

Our main goal is to further examine qualitative properties of the SD method relevant  with the stability of the method and answer questions of the following type: \textit{Is the SD method able to preserve the asymptotic stability of the underlying SDE?}


The answer of the question above is to the positive, and is given in our main result, Theorem \ref{asTSD:theorem:asNUMstability}.
In Section \ref{asTSD:sec:main} we give all the necessary information about the truncated version of the semi-discrete method; the way of construction of the numerical scheme and some useful properties, whereas Section \ref{asTSD:sec:stability} contains the main result with the proof. 
Section \ref{asTSD:sec:example} provides a numerical example. Motivated by the SDE appearing in the example, we also propose a different SD scheme, using the Lamperti transformation to the original SDE, which we call Lamperti semi-discrete (LSD). Numerical simulations support our theoretical findings. Finally, Section \ref{asTSD:sec:conclusion} contains concluding remarks.

\section{Setting and Assumptions}\label{asTSD:sec:main}

Throughout, let $T>0$ and $(\Omega, \bbf, \{\bbf_t\}_{0\leq t\leq T}, \bfP)$ be a complete probability space, meaning that the filtration $ \{\bbf_t\}_{0\leq t\leq T} $ satisfies the usual conditions, i.e. is right continuous and $\bbf_0$ includes all $\bfP$-null sets. Let $W_{t,\w}:[0,T]\times\W\rightarrow\bbR$ be a one-dimensional Wiener process adapted to the filtration $\{\bbf_t\}_{0\leq t\leq  T}.$  Consider SDE  (\ref{asTSD-eq:scalarSDEs}), which we rewrite here in its integral form
\beqq\label{asTSD-eq:general sde}
x_t=x_0 + \int_{0}^{t}a(s,x_s)ds + \int_{0}^{t}b(s,x_s)dW_s,\quad t\in [0,T],
\eeqq
which admits a unique strong solution. In particular, we assume the existence of a predictable stochastic process $x:[0,T]\times \W\rightarrow \bbR$ such that (\cite[Def. 2.1]{mao:2007}),
$$
\{a(t,x_t)\}\in\bbl^1([0,T];\bbR), \quad \{b(t,x_t)\}\in\bbl^2([0,T];\bbR)
$$
and
$$
\bfP\left[x_t=x_0 + \int_{0}^{t}a(s,x_s)ds + \int_{0}^{t}b(s,x_s)dW_s\right]=1, \quad \hbox{ for every } t\in[0,T].
$$

\bass\label{asTSD:assA}
Let $f(s,r,x,y), g(s,r,x,y):[0,T]^2\times\bbR^2\rightarrow\bbR$ be such that $f(s,s,x,x)=a(s,x), g(s,s,x,x)=b(s,x),$ where $f,g$ satisfy the following condition $(\phi\equiv f,g)$
$$
|\phi(s_1,r_1,x_1,y_1) - \phi(s_2,r_2,x_2,y_2)|\leq C_R \Big( |s_1-s_2| + |r_1-r_2| + |x_1-x_2| + |y_1-y_2| \Big)
$$
for any $R>0$ such that $|x_1|\vee|x_2|\vee|y_1|\vee|y_2|\leq R,$ where the quantity $C_R$ depends on $R$ and $x\vee y$ denotes the maximum of $x, y.$
\eass

Let us now recall the SD scheme. Consider the equidistant partition $0=t_0<t_1<...<t_N=T$ and $\D=T/N.$
We assume that for every $n\leq N-1,$ the following SDE
\beqq\label{asTSD-eq:SD scheme}
y_t=y_{t_n} + \int_{t_n}^{t} f(t_n, s, y_{t_n}, y_s)ds
+ \int_{t_n}^{t} g(t_n, s, y_{t_n}, y_s)dW_s,\quad t\in(t_n, t_{n+1}],
\eeqq
with $y_0=x_0$ a.s., has a unique strong solution.

In order to compare with the exact solution $x_t,$ which is a continuous time process, we consider the following interpolation process of the semi-discrete approximation, in a compact form,
\beqq\label{asTSD-eq:SD scheme compact}
y_t=y_0 + \int_{0}^{t}f(\hat{s}, s,y_{\hat{s}},y_s)ds +
\int_{0}^{t}g(\hat{s},s,y_{\hat{s}},y_s) dW_s,
\eeqq
where $\hat{s}=t_{n}$ when $s\in[t_n,t_{n+1}).$ Process (\ref{asTSD-eq:SD scheme compact}) has jumps at nodes $t_n.$ The first and third variable in $f, g$ denote the discretized part of the original SDE. We observe from (\ref{asTSD-eq:SD scheme compact}) that in order to solve for $y_t$, we have to solve an SDE and not an algebraic equation. The choice $f(s,r,x,y)=a(s,x)$ and $g(s,r,x,y)=b(s,x)$ reproduces the classical Euler scheme. 

In the case of superlinear coefficients the numerical scheme (\ref{asTSD-eq:SD scheme compact}) converges to the true solution $x_t$ of SDE (\ref{asTSD-eq:general sde}) and this is stated in the following, cf. \cite{halidias_stamatiou:2016},
\bth[Strong convergence]\label{asTSD-thm:strong_conv}
	Suppose Assumption \ref{asTSD:assA} holds and (\ref{asTSD-eq:SD scheme}) has  a unique strong solution for every $n\leq N-1,$ where $x_0\in \bbl^p(\Omega,\bbR).$ Let also
	$$
	\bfE(\sup_{0\leq t\leq T}|x_t|^p) \vee \bfE(\sup_{0\leq t\leq T}|y_t|^p)<A,
	$$
	for some $p>2$ and $A>0.$ Then the semi-discrete numerical scheme (\ref{asTSD-eq:SD scheme compact}) converges to the true solution of (\ref{asTSD-eq:general sde}) in the $\bbl^2$-sense, that is
	\beqq \label{asTSD-eq:strong_conv}
	\lim_{\D\rightarrow0}\bfE\sup_{0\leq t\leq T}|y_t-x_t|^2=0.
	\eeqq
\ethe

Relation (\ref{asTSD-eq:strong_conv}) does not reveal the order of convergence. We choose a strictly increasing function $\mu:\bbR_+\rightarrow \bbR_+$ such that 
for every $s,r\leq T$
\beqq \label{asTSD-eq:mu}
\sup_{|x|\leq u}\left(|f(s,r,x,y)| \vee |g(s,r,x,y)|\right)\leq \mu(u)(1 + |y|), \qquad u\geq1.
\eeqq

The inverse function of $\mu,$  denoted by $\mu^{-1},$ maps $[\mu(1),\infty)$ to $\bbR_+.$ Moreover, we choose a strictly decreasing function $h:(0,1]\rightarrow[\mu(1),\infty)$ and a constant $\hat{h}\geq 1\vee \mu(1)$ such that
\beqq \label{asTSD-eq:h}
\lim_{\D\rightarrow0}h(\D)=\infty \quad \hbox{and}\quad \D^{1/6}h(\D)\leq \hat{h} \quad \hbox{for every}\quad  \D\in(0,1].
\eeqq

Now, we are ready to define the truncated versions of $f, g.$ Let $\D\in(0,1]$ and $f_\D, g_\D$ defined by
\beqq \label{asTSD-eq:trunc}
\phi_\D(s,r,x,y):=\phi\left(s,r,(|x|\wedge\mu^{-1}(h(\D)))\frac{x}{|x|},y\right),
\eeqq
for $x,y\in\bbR$ where we set $x/|x|=0$ when $x=0.$

It follows that the truncated functions $f_\D, g_\D$ are bounded in the following way for a given step-size $0<\D\leq1,$
\beam \nonumber
|f_\D(s,r,x,y)| \vee |g_\D(s,r,x,y)|&\leq& \mu(\mu^{-1}(h(\D)))(1 + |y|)\\
\label{asTSD-eq:bounded_fg}&\leq &h(\D)(1 + |y|),
\eeam
for all $x,y\in\bbR.$

For the equidistant partition of $[0,T]$ with $\D<1$ consider now the following SDE
\beqq\label{asTSD-eq:SD scheme_trunc}
y_t^\D=y_{t_n}^\D + \int_{t_n}^{t} f_\D(t_n, s, y_{t_n}^\D, y_s^\D)ds
+ \int_{t_n}^{t} g_\D(t_n, s, y_{t_n}^\D, y_s^\D)dW_s,\quad t\in(t_n, t_{n+1}],
\eeqq
with $y_0=x_0$ a.s. We assume that (\ref{asTSD-eq:SD scheme_trunc}) admits a unique strong solution for every $n\leq N-1$ and rewrite it in compact form,
\beqq\label{asTSD-eq:SD_scheme_compact_trunc}
y_t^\D=y_0 + \int_{0}^{t}f_\D(\hat{s}, s,y_{\hat{s}}^\D,y_s^\D)ds +
\int_{0}^{t}g_\D(\hat{s},s,y_{\hat{s}}^\D,y_s^\D) dW_s.
\eeqq

\bass\label{asTSD:assB}
Let the truncated versions $f_\D(s,r,x,y), g_\D(s,r,x,y)$ of $f, g$  satisfy the following condition $(\phi_\D\equiv f_\D,g_\D)$
$$
|\phi_\D(s_1,r_1,x_1,y_1) - \phi_\D(s_2,r_2,x_2,y_2)|\leq h(\D) \Big( |s_1-s_2| + |r_1-r_2| + |x_1-x_2| + |y_1-y_2| \Big)
$$
for all $0<\D\leq 1$ and $x_1, x_2, y_1, y_2\in \bbR,$ where $h(\D)$ is as in (\ref{asTSD-eq:h}).
\eass

Let us also assume that the coefficients $a(t,x), b(t,x)$ of the original SDE satisfy the Khasminskii-type condition.

\bass\label{asTSD:assC}
We assume the existence of constants $p\geq2$ and $C_K>0$ such that $x_0\in \bbl^p(\Omega,\bbR)$  and
$$
xa(t,x) + \frac{p-1}{2}b(t,x)^2\leq C_K(1 + |x|^2)
$$
for all $(t,x)\in[0,T]\times\bbR$.
\eass

A well-known result follows (see e.g. \cite{mao:2007}) when the SDE (\ref{asTSD-eq:general sde}) satisfies the local Lipschitz condition plus the Khasminskii-type condition.

\ble
	Under Assumptions \ref{asTSD:assA} (for the coefficients $a(t,x), b(t,x)$) and \ref{asTSD:assC} the SDE (\ref{asTSD-eq:general sde}) has a unique global solution and for all $T>0,$ there exists a constant $A>0$ such that
	$$
	\sup_{0\leq t\leq T}\bfE |x_t|^p<A.
	$$
\ele

\bth[Order of strong
convergence]\label{asTSD:theorem:StrongConvergenceOrder} Suppose
Assumption \ref{asTSD:assB} and Assumption \ref{asTSD:assC} hold and
(\ref{asTSD-eq:SD scheme_trunc}) has a unique strong solution for
every $n\leq N-1,$ where $x_0\in \bbl^p(\Omega,\bbR)$ for some
$p\geq 14+2\ga.$ Let $\ep\in(0,1/3)$ and define for $\ga>0$
$$
\mu(u) = \ov{C}u^{1+\ga}, \quad u\geq0 \quad \hbox { and } \quad
h(\D)=\ov{C} + \sqrt{\ln \D^{-\ep}}, \quad\D\in(0,1].
$$
where $\D\leq1$ and $\hat{h}$ are such that (\ref{asTSD-eq:h})
holds. Then the semi-discrete numerical scheme (\ref{asTSD-eq:SD_scheme_compact_trunc}) converges to the true solution of (\ref{asTSD-eq:general sde}) in the $\bbl^2$-sense with order
arbitrarily close to $1/2,$ that is 
\beqq
\label{asTSD-eq:strong_conv_order} \bfE\sup_{0\leq t\leq
	T}|y^\D_t-x_t|^2\leq C\D^{1-\ep}. \eeqq 
\ethe

\section{Asymptotic Stability}\label{asTSD:sec:stability}
%

Now we are ready to study the ability of the truncated SD method to preserve the asymptotic stability of (\ref{asTSD-eq:general sde}). For that reason we also assume that $a(0,0) = 0$ and $b(0,0)=0.$ Moreover, to guarantee the asymptotic stability of (\ref{asTSD-eq:general sde}) we use an assumption similar to \cite[Assumption 5.1]{hu_li_mao:2018}. 
\bass\label{asTSD:assD}
We assume the existence of a continuous non-decreasing function $\ka:\bbR_+\mapsto\bbR_+$ with $\ka(0)=0$ and $\ka(u)>0$ for all $u>0$ such that
\beqq\label{asTSD-eq:asstabcond}
2xa(s,x) + |b(x)|^2\leq - \ka(|x|),
\eeqq
for all $x\in\bbR$ and $s\in[0,T].$
\eass
Now, we state a result without proof concerning the asymptotic stability of (\ref{asTSD-eq:general sde}), see also \cite[Theorem  5.2]{hu_li_mao:2018} where autonomous coefficients are assumed.

\bth[asymptotic stability of underlying process]\label{asTSD:theorem:assstability}
Let Assumption \ref{asTSD:assD} hold. Then the solution process of SDE (\ref{asTSD-eq:general sde}) is asymptotically stable, that is
\beqq\label{asTSD-eq:asstab}
\lim_{\tto}x_t =0  \hbox{ a.s.}
\eeqq
for any $x_0\in\bbR.$ 
\ethe

Recall equation (\ref{asTSD-eq:SD scheme_trunc}) which defines the truncated SD numerical scheme. We rewrite our proposed scheme, that is the solution of (\ref{asTSD-eq:SD scheme_trunc}) at the discrete points $0,t_1,\ldots,t_{n+1},$ in the following way
\beqq\label{asTSD-eq:SDrepresentation}
y_{n+1}^\D=\phi^\D(y_{n}^\D,t_n,\D,\D W_n), 
\eeqq
where $\D W_n$ are the Wiener increments, $\D=t_{n+1}-t_n$ is the step-size and $y_n$ stands for $y_{t_n}.$ We assume the following decomposition of $\phi^\D(y_{n}^\D,t_n,\D,\D W_n)$ for the above representation (\ref{asTSD-eq:SDrepresentation}),
\beqq \label{asTSD-eq:NUMrepresentation}
\phi^\D(y_{n}^\D,t_n,\D,\D W_n)^2 = (y_{n}^\D)^2 +  \phi_1^\D(y_{n}^\D,t_n,\D) + \phi_2^\D(y_{n}^\D,t_n,\D,\D W_n),
\eeqq
where $\bfE(\phi_2^\D(y_{n}^\D,t_n,\D,\D W_n)|\bbf_{t_n})=0.$
The following theorem shows that the truncated SD method is able to preserve the asymptotic stability property of the underlying SDE.

\bth[asymptotic numeric stability]\label{asTSD:theorem:asNUMstability}
Let the auxiliary function $\phi_1^\D$ from (\ref{asTSD-eq:NUMrepresentation}) satisfy 
\beqq \label{asTSD-eq:asstab_num_cond_phi}
\phi_1^\D(y_{n}^\D,t_n,\D) \leq -\ka_1\left(\left|(|y_{n}^\D|\wedge\mu^{-1}(h(\D)))\frac{y_{n}^\D}{|y_{n}^\D|}\right|\right),
\eeqq
for any $0<\D\leq\D^*,$ where $\ka_1$ has the same properties as $\ka$ in (\ref{asTSD-eq:asstabcond}) with $\ka_1\leq\ka$. Let also Assumption \ref{asTSD:assD} hold. 

Then the solution of the truncated SD method (\ref{asTSD-eq:SDrepresentation}) is numerically asymptotically stable, that is
\beqq \label{asTSD-eq:asNUMstab}
\lim_{\nto}y_n^{\D} =0  \hbox{ a.s.}
\eeqq
for all $x_0\in\bbR$ and $0<\D\leq\D^*$.
\ethe

\bpf[Proof of Theorem \ref{asTSD:theorem:asNUMstability}]
Let us first fix a $\D\in(0,\D^*].$ Denote $$\pi_\D(x) := \left(|x|\wedge\mu^{-1}(h(\D))\right)\frac{x}{|x|}.$$ Then combining (\ref{asTSD-eq:SDrepresentation}), (\ref{asTSD-eq:NUMrepresentation}) and (\ref{asTSD-eq:asstab_num_cond_phi}) we get
\beao
(y_{n+1}^\D)^2 &\leq& (y_{n}^\D)^2 -\ka_1\left(\left|\pi_\D(y_{n}^\D)\right|\right) + \phi_2^\D(y_{n}^\D,t_n,\D,\D W_n)\\
&\leq& (x_{0})^2 -\sum_{j=0}^n\ka_1\left(\left|\pi_\D(y_{j}^\D)\right|\right) + M_n,
\eeao
where $M_n:=\sum_{j=0}^n\phi_2^\D(y_{j}^\D,t_j,\D,\D W_j).$ Recalling that $\bfE(\phi_2^\D(y_{n}^\D,t_n,\D,\D W_n)|\bbf_{t_n})=0$ implies that $M_n, n=0,1,\ldots,$ is a martingale. Application of the nonnegative semi-martingale convergence theorem, c.f. \cite[Theorem 7, p.139]{lipster_shiryayev:1989}, implies
$$\sum_{j=0}^\infty\ka_1\left(\left|\pi_\D(y_{j}^\D)\right|\right)<\infty \hbox{ a.s.}$$
which in turn 
$$\lim_{j\rightarrow\infty}\ka_1\left(\left|\pi_\D(y_{j}^\D)\right|\right)=0\hbox{ a.s.}$$
By the property of the function $\ka_1$ we get that 
$$\lim_{j\rightarrow\infty}\left(|y_{j}^\D|\wedge\mu^{-1}(h(\D))\right)\frac{y_{j}^\D}{|y_{j}^\D|}=0\hbox{ a.s.}$$
Assertion (\ref{asTSD-eq:asNUMstab}) follows.
\epf

\section{Example}\label{asTSD:sec:example}

We will use the numerical example of \cite[Example 5.4]{hu_li_mao:2018}, that is we consider an autonomous SDE of the form (\ref{asTSD-eq:general sde}) with $a(x)=-10x^3$ and $b(x)=x^2,$ with initial condition  $x_0\in\bbR$, that is,
\beqq  \label{asTSD-eq:exampleSDE}
x_t =x_0 -10\int_0^t x_s^3 ds + \int_0^t x_s^2 dW_s, \qquad t\geq0.
\eeqq

Using standard arguments one may show that the solution process of SDE (\ref{asTSD-eq:exampleSDE}) is positive, see Appendix \ref{asTSD-ap:positivity_of_SDE}. 
Assumption \ref{asTSD:assD} holds with $\ka(u) = 19u^4$ therefore by Theorem \ref{asTSD:theorem:assstability} SDE (\ref{asTSD-eq:exampleSDE}) is almost surely asymptotically stable. The classical Euler Maruyama method is not able to reproduce this asymptotic stability, see \cite[Appendix]{hu_li_mao:2018}. In the following we show that  the truncated SD method can reproduce this asymptotic stability. Since, in the construction of the semi discrete method the way of discretizing is not unique (but rather indicated by the equation itself) we will try two versions of the SD method by freezing different parts of the diffusion coefficient.
We first choose the auxiliary functions $f, g_1$ and $g_2$ in the following way
$$
f(s,r,x,y) = -10x^2y, \qquad g_1(s,r,x,y) = x,\qquad g_2(s,r,x,y) = x^2,
$$
thus (\ref{asTSD-eq:SD scheme}) becomes
\beqq\label{asTSD-eq:SD schemeExample}
y_t=y_{t_n} - 10y_{t_n}^2\int_{t_n}^{t} y_sds + y_{t_n}^2\int_{t_n}^{t}dW_s, \quad t\in(t_n, t_{n+1}],
\eeqq
and
\beqq\label{asTSD-eq:SD scheme2Example}
\hat{y}_t = \hat{y}_{t_n} - 10\hat{y}_{t_n}^2\int_{t_n}^{t} \hat{y}_sds + \hat{y}_{t_n}\int_{t_n}^{t}\hat{y}_sdW_s, \quad t\in(t_n, t_{n+1}],
\eeqq
respectively, with $y_0=\hat{y}_0 = x_0$ a.s. SDEs (\ref{asTSD-eq:SD schemeExample}) and (\ref{asTSD-eq:SD scheme2Example}) are linear equations ( (\ref{asTSD-eq:SD schemeExample}) is linear in the narrow sense and is known as Langevin equation) with variable coefficients which admit a unique strong solution, c.f. \cite[Chapter 4.4]{kloeden_platen:1995} and Appendix \ref{asTSD-ap:solution_linear_SDEs}. In particular,
\beqq  \label{asTSD-eq:exampleSD}
y_{n+1} =e^{-10y_{n}^2\D}\left(y_n + y_n^2\int_{t_n}^{t_{n+1}}e^{10y_{n}^2(s-t_n)}dW_s\right), \quad n\in\bbN,
\eeqq
and
\beqq  \label{asTSD-eq:exampleSD2}
\hat{y}_{n+1} =\hat{y}_n\exp\left\{-\frac{21}{2}\hat{y}_{n}^2\D + \hat{y}_n\D W_n\right\}, \quad n\in\bbN.
\eeqq

Note that (\ref{asTSD-eq:mu}) holds with $\mu(u) =10|u|^2$ since
$$
\sup_{|x|\leq u}\left(|-10x^2y| \vee |x|\vee |x^2|\right)\leq 10|u|^2(1 + |y|), \qquad u\geq1.
$$
Therefore, in the notation of Theorem \ref{asTSD:theorem:StrongConvergenceOrder}, $\ga=1$ and $\ov{C}=10.$ Finally, $h(\D) = 10 + \sqrt{\ln \D^{-\ep_1}}$ for any $\D\in(0,1].$ Clearly $h(1)\geq \mu(1)$ and
$$
\D^{1/6}h(\D) \leq 10\D^{1/6} + \sqrt{\D^{1/3}\ln \D^{-\ep_1}}\leq 11,
$$
for any $\D\in(0,1]$ and $0<\ep_1\leq 1/3.$ Therefore we take $\hat{h}=11.$
The truncated versions of the semi-discrete method (TSD) read,
\beqq  \label{asTSD-eq:exampleSDtrunc}
y_{n+1}^\D =e^{-10\pi^2_\D(y_{n}^\D)\D}\left(y_{n}^\D + \pi^2_\D(y_{n}^\D)\int_{t_n}^{t_{n+1}}e^{10\pi^2_\D(y_{n}^\D)(s-t_n)}dW_s\right)
\eeqq
and
\beqq  \label{asTSD-eq:exampleSD2trunc}
\hat{y}_{n+1}^\D =\hat{y}_{n}^\D\exp\left\{-\frac{21}{2}\pi^2_\D(\hat{y}_{n}^\D)\D + \pi_\D(\hat{y}_{n}^\D)\D W_n\right\}
\eeqq
for $n\in\bbN,$ where 
$$\pi_\D(x)=\left(|x|\wedge \sqrt{\frac{h(\D)}{10}}\right)\frac{x}{|x|}$$ and therefore
$$\pi^2_\D(x)=|x|^2\wedge\frac{h(\D)}{10}.$$

\subsection{Asymptotic stability of truncated Semi-Discrete method}\label{asTSD:ssec:tsd}

Now, we compute $(y_{n+1}^\D)^2$ taking the square of (\ref{asTSD-eq:exampleSDtrunc})  and making some rearrangements to show that it admits representation (\ref{asTSD-eq:NUMrepresentation}). 
\beao
&&(y_{n+1}^\D)^2 = (y_{n}^\D)^2e^{-20\pi^2_\D(y_{n}^\D)\D} + 2y_{n}^\D e^{-10\pi^2_\D(y_{n}^\D)t_{n+1}}\pi^2_\D(y_{n}^\D)\int_{t_n}^{t_{n+1}}e^{10\pi^2_\D(y_{n}^\D)s}dW_s\\
&&+ \pi^4_\D(y_{n}^\D)e^{-20\pi^2_\D(y_{n}^\D)t_{n+1}}\left(\int_{t_n}^{t_{n+1}}e^{10\pi^2_\D(y_{n}^\D)s}dW_s\right)^2\\
&=& (y_{n}^\D)^2 - (1-e^{-20\pi^2_\D(y_{n}^\D)t_{n+1}})(y_{n}^\D)^2 + \frac{1}{20}\pi^2_\D(y_{n}^\D)e^{-20\pi^2_\D(y_{n}^\D)t_{n+1}}(e^{20\pi^2_\D(y_{n}^\D)t_{n+1}}-e^{20\pi^2_\D(y_{n}^\D)t_{n}}) \\
&& + 2y_{n}^\D e^{-10\pi^2_\D(y_{n}^\D)t_{n+1}}\pi^2_\D(y_{n}^\D)\int_{t_n}^{t_{n+1}}e^{10\pi^2_\D(y_{n}^\D)s}d W_s\\
&& + \pi^4_\D(y_{n}^\D)e^{-20\pi^2_\D(y_{n}^\D)t_{n+1}}\left(\int_{t_n}^{t_{n+1}}e^{10\pi^2_\D(y_{n}^\D)s}d W_s\right)^2\\
&& - \frac{1}{20}\pi^2_\D(y_{n}^\D)e^{-20\pi^2_\D(y_{n}^\D)t_{n+1}}(e^{20\pi^2_\D(y_{n}^\D)t_{n+1}}-e^{20\pi^2_\D(y_{n}^\D)t_{n}}).
\eeao

Denote $I := I(y_{n}^\D,t_n,\D,\D W_n) = \int_{t_n}^{t_{n+1}}e^{10\pi^2_\D(y_{n}^\D)s}d W_s$ and set 
\beao
\phi_2^\D(y_{n}^\D,t_n,\D,\D W_n)&:=&2y_{n}^\D e^{-10\pi^2_\D(y_{n}^\D)t_{n+1}}\pi^2_\D(y_{n}^\D)I + \pi^4_\D(y_{n}^\D)e^{-20\pi^2_\D(y_{n}^\D)t_{n+1}}I^2\\
&& -\frac{1}{20}\pi^2_\D(y_{n}^\D)e^{-20\pi^2_\D(y_{n}^\D)t_{n+1}}(e^{20\pi^2_\D(y_{n}^\D)t_{n+1}}-e^{20\pi^2_\D(y_{n}^\D)t_{n}})
\eeao 
to see that $\bfE(\phi_2^\D(y_{n}^\D,t_n,\D,\D W_n)|\bbf_{t_n})=0.$
Moreover 
\beao
\phi_1^\D(y_{n}^\D,t_n,\D)&:=&-(1-e^{-20\pi^2_\D(y_{n}^\D)\D})(y_{n}^\D)^2 + \frac{1}{20}\pi^2_\D(y_{n}^\D)(1-e^{-20\pi^2_\D(y_{n}^\D)\D})\\
&\leq&-\frac{19}{20}\left(1-e^{-20\pi^2_\D(y_{n}^\D)\D}\right)\pi^2_\D(y_{n}^\D),  
\eeao
implying that we may choose $\ka_1$ in the following way
$$
\ka_1(u):=-\frac{19}{20}(1-e^{-20u^2\D})u^2,
$$
so that condition (\ref{asTSD-eq:asstab_num_cond_phi})  holds and therefore Theorem \ref{asTSD:theorem:asNUMstability}  applies. Note that $\ka_1(0)=0$ and $\ka_1(u)>0$ for any $\D>0.$ 
We conclude that the truncated SD scheme (\ref{asTSD-eq:exampleSDtrunc}) preserves the asymptotic stability perfectly in the sense that $\lim_{\nto}y_n^{\D} = 0$ a.s. for any $0<\D\leq1.$

\subsection{Asymptotic stability of exponential truncated Semi-Discrete method}\label{asTSD:ssec:exptsd}

We examine $(\hat{y}_{n+1}^\D)^2.$ We take the square of (\ref{asTSD-eq:exampleSD2trunc}) and get that 
\beao
&&(\hat{y}_{n+1}^\D)^2 =(\hat{y}_{n}^\D)^2e^{-21\pi^2_\D(\hat{y}_{n}^\D)\D + 2\pi_\D(\hat{y}_{n}^\D)\D W_n}\\
&=& (\hat{y}_{n}^\D)^2 - (1-e^{-19\pi^2_\D(\hat{y}_{n}^\D)\D})(\hat{y}_{n}^\D)^2 + (\hat{y}_{n}^\D)^2e^{-19\pi^2_\D(\hat{y}_{n}^\D)\D}(1-e^{-2\pi^2_\D(\hat{y}_{n}^\D)\D + 2\pi_\D(\hat{y}_{n}^\D)\D W_n}). 
\eeao
Set the last term of the above equality to $\phi_2^\D,$ that is
$$
\phi_2^\D(\hat{y}_{n}^\D,t_n,\D,\D W_n):=(\hat{y}_{n}^\D)^2e^{-19\pi^2_\D(\hat{y}_{n}^\D)\D}(1-e^{-2\pi^2_\D(\hat{y}_{n}^\D)\D + 2\pi_\D(\hat{y}_{n}^\D)\D W_n})$$
to see that $\bfE(\phi_2^\D(y_{n}^\D,t_n,\D,\D W_n)|\bbf_{t_n})=0$ since $\bbE_n:=e^{-2\pi^2_\D(\hat{y}_{n}^\D)\D + 2\pi_\D(\hat{y}_{n}^\D)\D W_n}$ is an exponential martingale. 

Moreover 
\beao
\phi_1^\D(\hat{y}_{n}^\D,t_n,\D)&:=&-(1-e^{-19\pi^2_\D(\hat{y}_{n}^\D)\D})(\hat{y}_{n}^\D)^2\\
&\leq&-\left(1-e^{-19\pi^2_\D(\hat{y}_{n}^\D)\D}\right)\pi^2_\D(\hat{y}_{n}^\D),  
\eeao
implying that we may choose $\ka_1$ in the following way
$$
\ka_1(u):=-(1-e^{-19u^2\D})u^2,
$$
so that once more condition (\ref{asTSD-eq:asstab_num_cond_phi})  holds and consequently Theorem \ref{asTSD:theorem:asNUMstability}  applies. We conclude that the truncated exponential SD scheme (\ref{asTSD-eq:exampleSD2trunc}) preserves the asymptotic stability perfectly in the sense that $\lim_{\nto}\hat{y}_n^{\D} = 0$ a.s. for any $0<\D\leq1.$

\subsection{Semi-Discrete method and Lampreti transformation}\label{asTSD:ssec:Lamperti}

Instead of approximating directly (\ref{asTSD-eq:exampleSDE}) we first study a transformation of it, which produces a new SDE with constant diffusion coefficient; in other words we use the Lamperti transformation of (\ref{asTSD-eq:exampleSDE}). In particular, consider $z = -1/x.$ 
 The It\^o formula implies the following dynamics for $(z_t),$ see Appendix \ref{asTSD-ap:Lamperti_tranformation},
 
\beqq  \label{asTSD-eq:exampleSDELamperti}
z_t =z_0 + 11\int_0^t (z_s)^{-1} ds + \int_0^t  dW_s, \qquad t\geq0.
\eeqq 

Let $t\in(t_n, t_{n+1}]$ and 
\beqq\label{asTSD-eq:SD schemeExampleLT}
\wt{y}_t = \D W_n + \wt{y}_{t_n} + 11\int_{t_n}^{t} (\wt{y}_s)^{-1}ds,
\eeqq
with $\wt{y}_{0}=z_0.$ (\ref{asTSD-eq:SD schemeExampleLT}) is a Bernoulli type equation with solution satisfying
\beqq\label{asTSD-eq:SD schemeExampleLTsol}
(\wt{y}_t)^2 = (\D W_n + \wt{y}_{t_n})^2 + 22(t-t_n). 
\eeqq

Recall that when $x_0>0,$ the solution process  $x_t>0$ a.s. which implies $z_t<0$ a.s. which in turn suggests that we take the negative root of  (\ref{asTSD-eq:SD schemeExampleLTsol}) as the solution
Therefore we propose the following semi-discrete method for the approximation of (\ref{asTSD-eq:exampleSDELamperti}),  
\beqq\label{asTSD-eq:SD schemeExampleLT_transf}
\wt{y}_{t_{n+1}} = -\sqrt{(\D W_n + \wt{y}_{t_n})^2 + 22\D},
\eeqq
which suggests the Lamperti semi-discrete method $(\wt{z}_n)_{n\in\bbN}$ for the approximation of (\ref{asTSD-eq:exampleSDE})
\beqq\label{asTSD-eq:SD schemeExampleLToriginal}
\wt{z}_{t_{n+1}} = \frac{1}{\sqrt{(\D W_n + \wt{y}_{t_n})^2 + 22\D}}.
\eeqq

\subsection{Simulation Paths}\label{asTSD:ssecnumerics}

\begin{figure}[ht]
	\centering
	\begin{subfigure}{.45\textwidth}
		\includegraphics[width=1\textwidth]{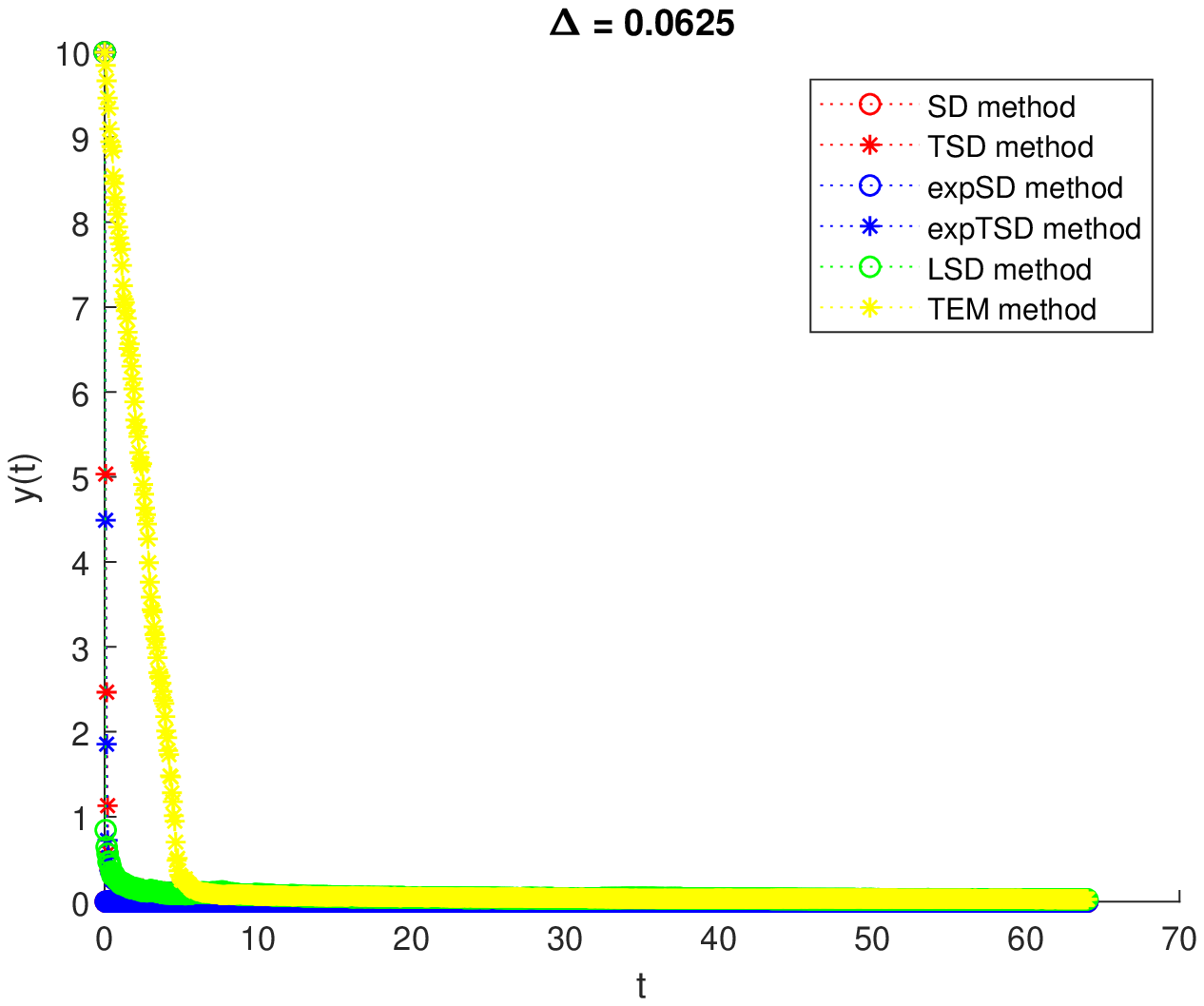}
		\caption{Trajectory for (\ref{asTSD-eq:exampleSD}) - (\ref{asTSD-eq:exampleSD2trunc}), (\ref{asTSD-eq:SD schemeExampleLToriginal}) and (\ref{asTSD-eq:exampleSDtruncEM}).}
	\end{subfigure}
	\begin{subfigure}{.45\textwidth}
		\includegraphics[width=1\textwidth]{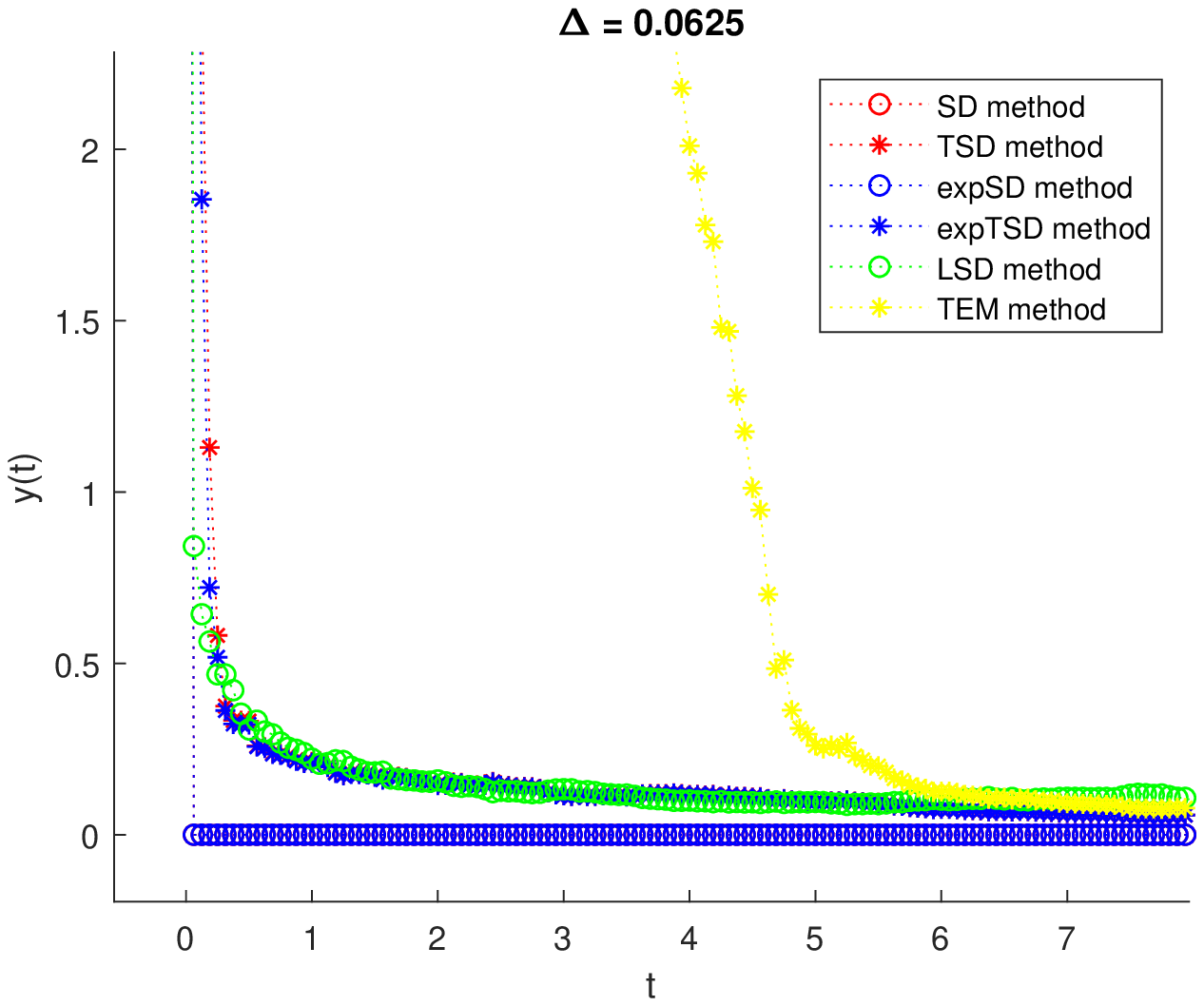}
		\caption{Zoom of Figure \ref{asTSD-fig:TSD_TEM}(A).}
	\end{subfigure}
	\caption{Trajectories of (\ref{asTSD-eq:exampleSD}) -(\ref{asTSD-eq:exampleSDtruncEM}), (\ref{asTSD-eq:SD schemeExampleLToriginal}) and (\ref{asTSD-eq:exampleSDtruncEM}) for the approximation of (\ref{asTSD-eq:exampleSDE}) with $\D<0.095$.}\label{asTSD-fig:TSD_TEM}
\end{figure}

\begin{figure}[ht]
	\centering
	\begin{subfigure}{.45\textwidth}
		\includegraphics[width=1\textwidth]{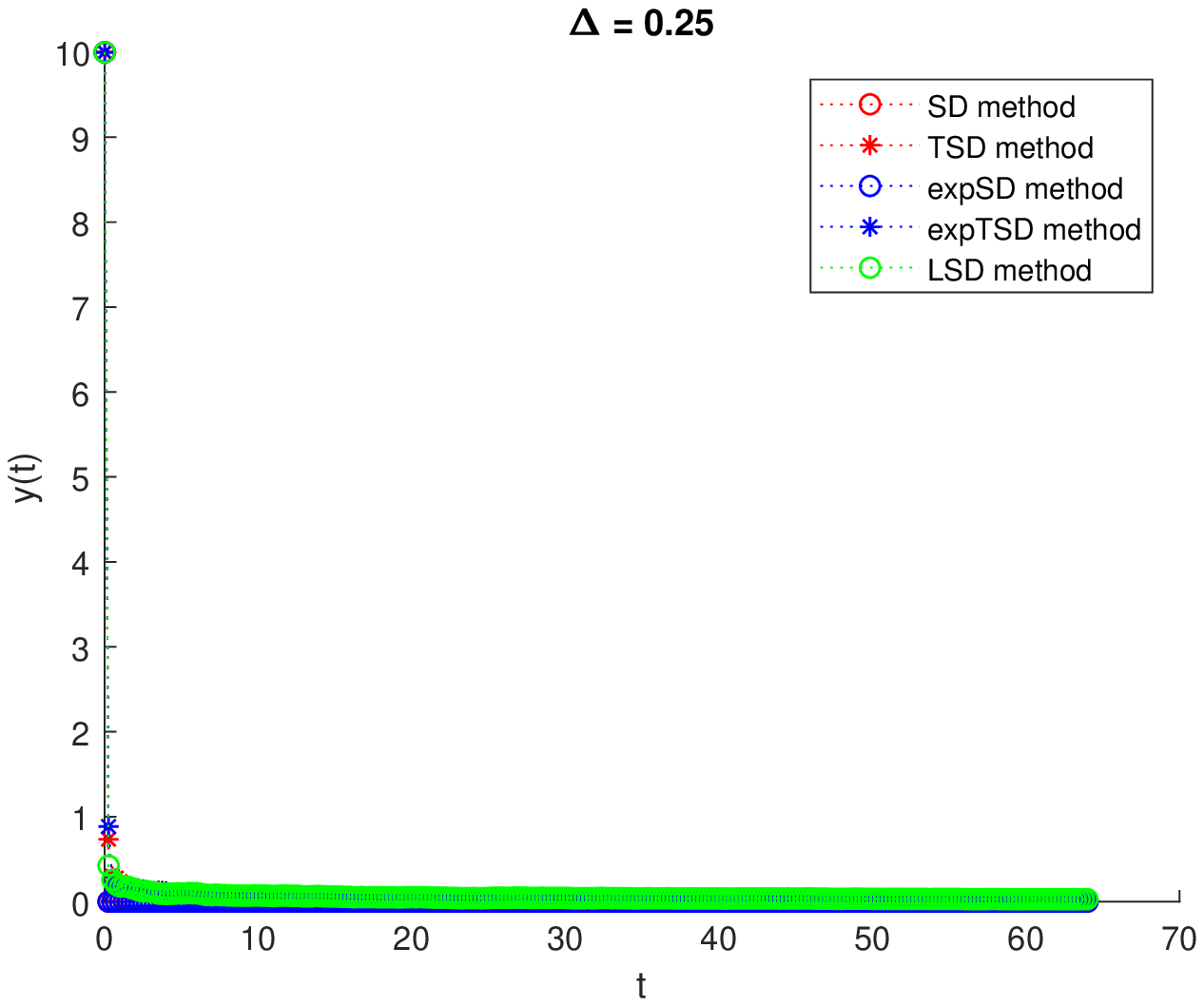}
		\caption{Trajectory for (\ref{asTSD-eq:exampleSD}) -(\ref{asTSD-eq:exampleSD2trunc}) and (\ref{asTSD-eq:SD schemeExampleLToriginal}).}
	\end{subfigure}
	\begin{subfigure}{.45\textwidth}
		\includegraphics[width=1\textwidth]{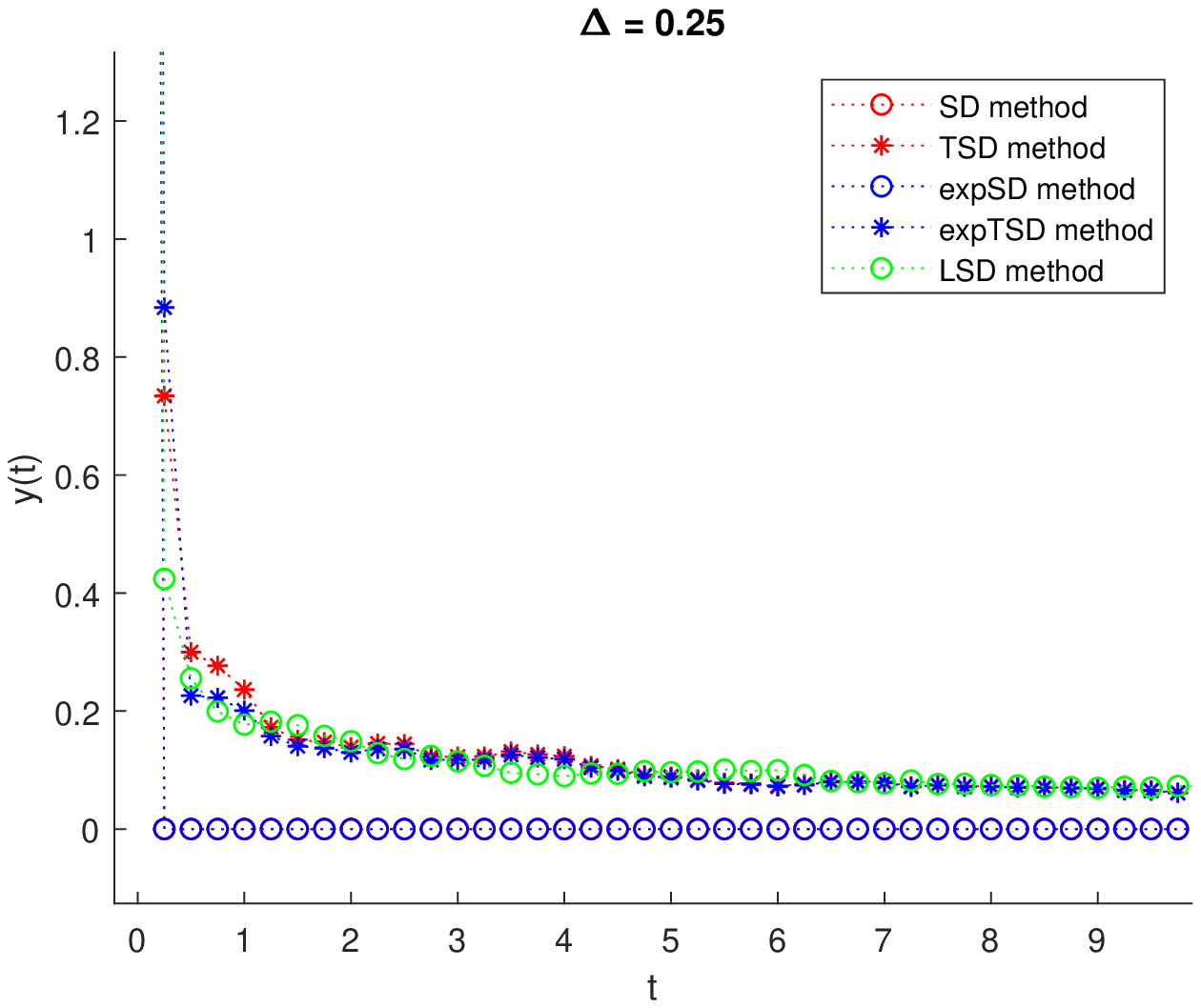}
		\caption{Zoom of Figure \ref{asTSD-fig:TSD}(A).}
	\end{subfigure}
	\caption{Trajectories of (\ref{asTSD-eq:exampleSD}) -(\ref{asTSD-eq:exampleSD2trunc}) and (\ref{asTSD-eq:SD schemeExampleLToriginal}) for the approximation of (\ref{asTSD-eq:exampleSDE}) with $\D=0.25$.}\label{asTSD-fig:TSD}
\end{figure}

\begin{figure}[ht]
	\centering
	\begin{subfigure}{.45\textwidth}
		\includegraphics[width=1\textwidth]{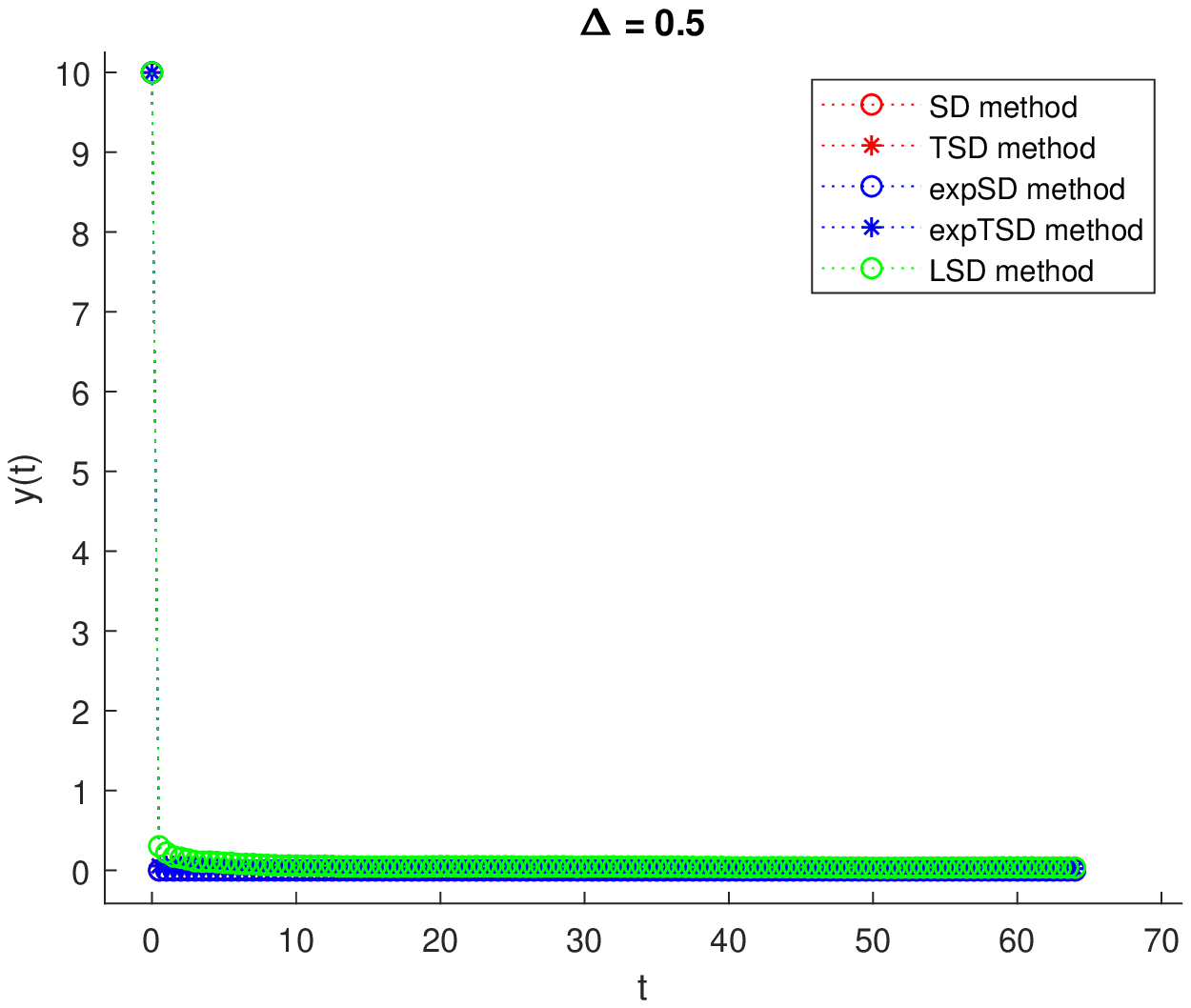}\label{asTSD-fig:TSD_TEM1}
		\caption{Trajectory for (\ref{asTSD-eq:exampleSD}) -(\ref{asTSD-eq:exampleSD2trunc}) and (\ref{asTSD-eq:SD schemeExampleLToriginal}).}
	\end{subfigure}
	\begin{subfigure}{.45\textwidth}
		\includegraphics[width=1\textwidth]{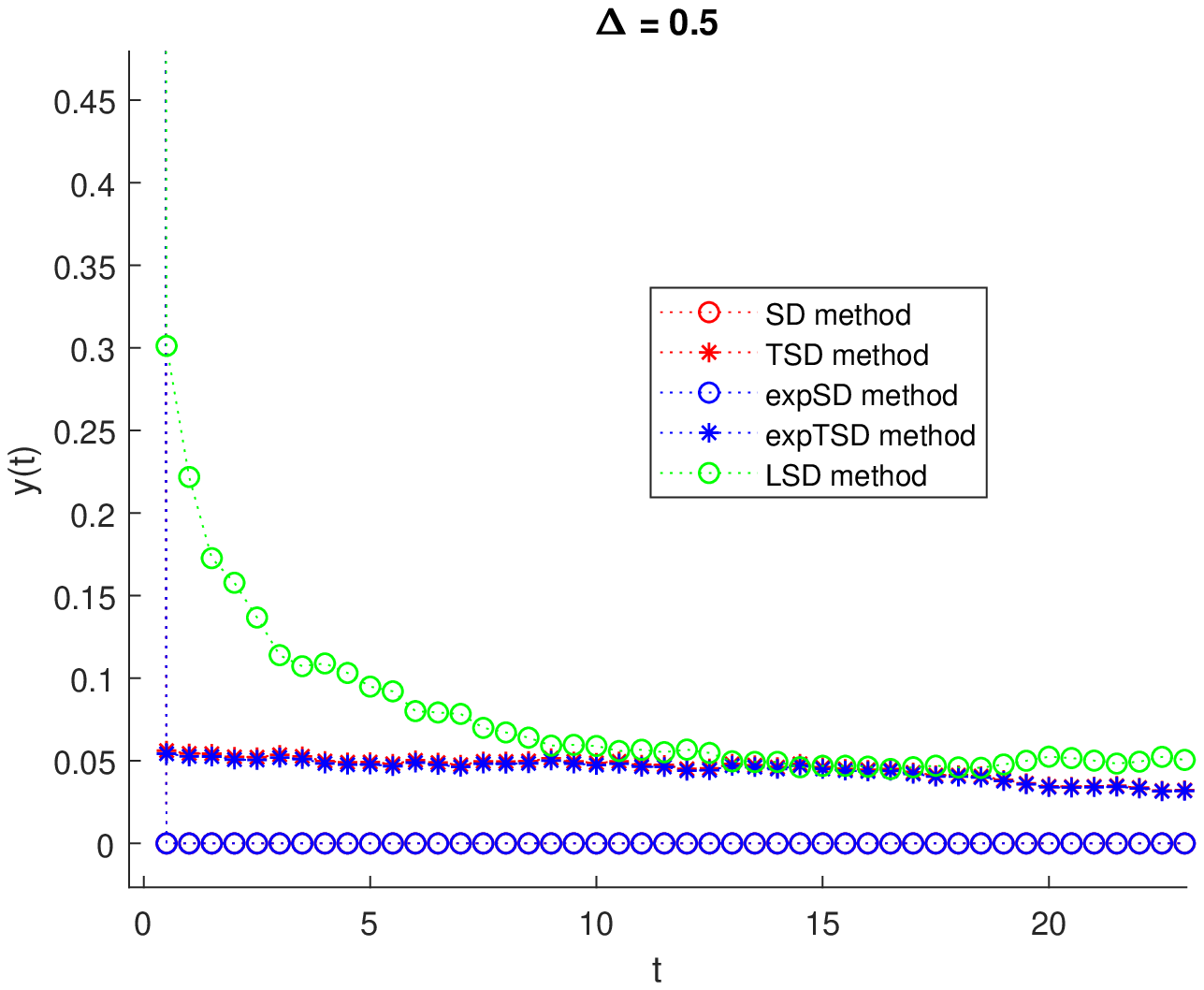}
		\caption{Zoom of Figure \ref{asTSD-fig:TSD2}(A).}
	\end{subfigure}
	\caption{Trajectories of (\ref{asTSD-eq:exampleSD}) - (\ref{asTSD-eq:exampleSD2trunc}) and (\ref{asTSD-eq:SD schemeExampleLToriginal}) for the approximation of (\ref{asTSD-eq:exampleSDE}) with $\D=0.5$.}\label{asTSD-fig:TSD2}
\end{figure}

We present simulations for the numerical approximation of (\ref{asTSD-eq:exampleSDE}) with  $x_0=10$ and compare with the truncated Euler Maruyama method (TEM), which reads
\beqq  \label{asTSD-eq:exampleSDtruncEM}
y_{n+1}^{TEM} =y_n - 10\left(|y_{n}|\wedge \bar{\mu}^{-1}(\bar{h}(\D))\frac{y_n}{|y_n|}\right)^3\D + \left(|y_{n}|\wedge \bar{\mu}^{-1}(\bar{h}(\D))\right)^2 \D W_n,
\eeqq
for $n\in\bbN,$ where $\bar{h}(\D)=\D^{-1/4}, \bar{\mu}(u) = 10u^3.$ According to the results in \cite{hu_li_mao:2018} it is shown that method (\ref{asTSD-eq:exampleSDtruncEM}) is asymptotically stable for any $\D\leq0.095,$ therefore for such small step sizes we compare all the methods presented here and for bigger $\D$ only the SD schemes (\ref{asTSD-eq:exampleSD}), (\ref{asTSD-eq:exampleSD2}), (\ref{asTSD-eq:exampleSDtrunc}) and (\ref{asTSD-eq:exampleSD2trunc}). We also present the Lamperti semi-discrete scheme (LSD) (\ref{asTSD-eq:SD schemeExampleLToriginal}). Moreover, the TEM method does not preserve positivity. Figures \ref{asTSD-fig:TSD_TEM}, \ref{asTSD-fig:TSD} and \ref{asTSD-fig:TSD2} shows sample simulations paths of TSD and TEM respectively for various stepsizes. Note that the truncated TSD, exponential truncated expTSD and the Lamperti LSD works for all  $\D< 1.$  

In the numerical simulation of the stochastic integral of the (truncated) TSD methods (\ref{asTSD-eq:exampleSD}) and (\ref{asTSD-eq:exampleSDtrunc}) we used the approximation
$\int_{t_n}^{t_{n+1}}e^{10\pi^2_\D(y_{n}^\D)s}dW_s\approx e^{10\pi^2_\D(y_{n}^\D)t_n}\D W_n,$ that is we calculated the integrand in the lower limit of integration. The above equality is of the order $\D^{r},$ with $0<r<1/2,$ see Appendix \ref{asTSD-ap:integralapprox}.

\begin{figure}[ht]
	\centering
		\includegraphics[width=1\textwidth]{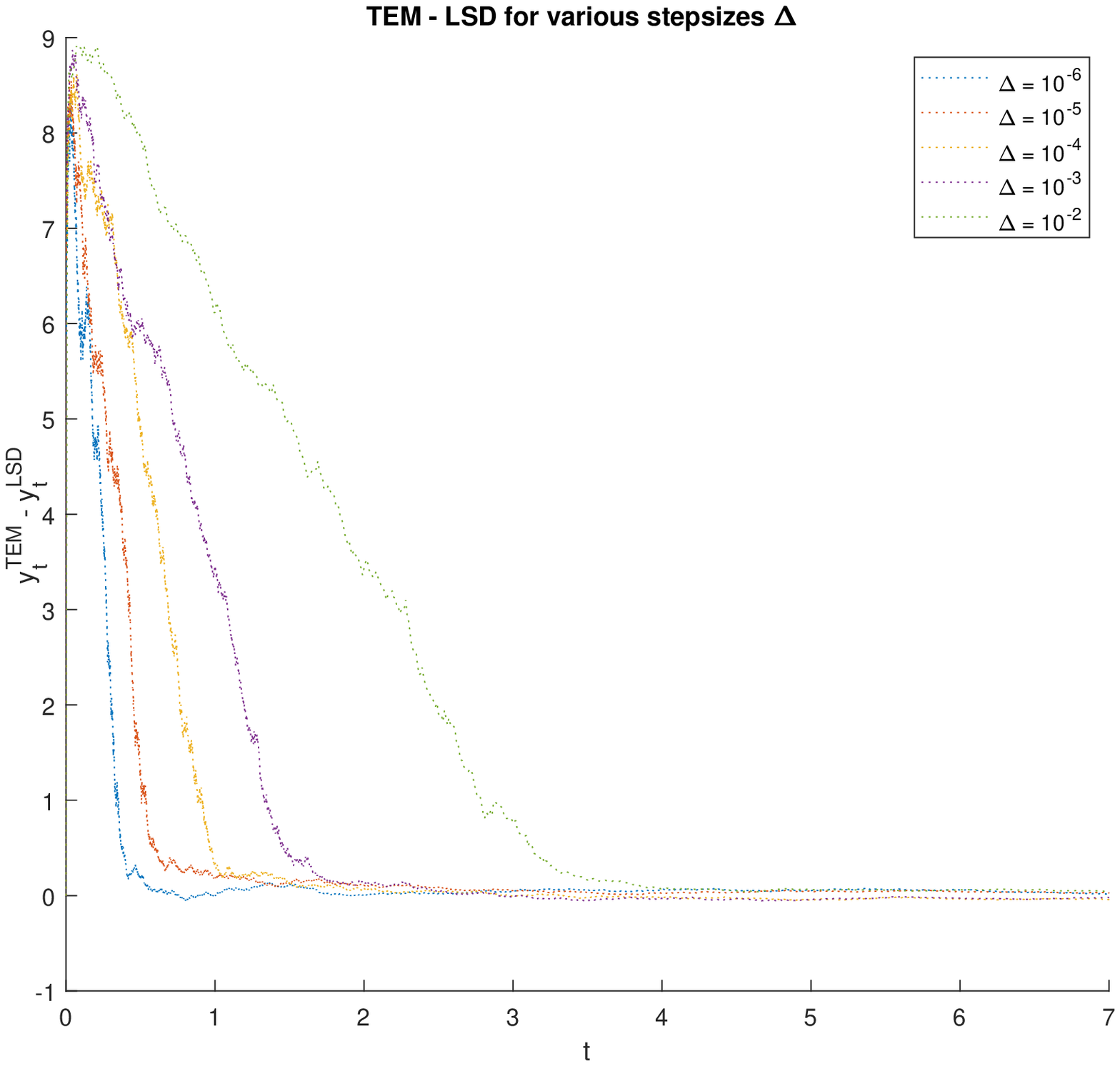}
		\caption{Difference of  (\ref{asTSD-eq:exampleSDtruncEM}) - (\ref{asTSD-eq:SD schemeExampleLToriginal}) for the approximation of (\ref{asTSD-eq:exampleSDE}) with various step sizes.}\label{asTSD-fig:TEMminTSD}
\end{figure}

We also present in Figure \ref{asTSD-fig:TEMminTSD} the difference between the Lamperti semi-discrete and the truncated Euler-Maryauma scheme, $y_{TEM}-y_{TSD}$ for small enough $\D$ such that $TSD$ works.

\section{Conclusion and Future Work}\label{asTSD:sec:conclusion}
In this paper we study the asymptotic stability of the semi-discrete (SD) numerical method for the approximation of stochastic differential equations. Recently, we examined the order of $\bbl^2$-convergence of the truncated SD method and showed that it can be arbitrarily close to $1/2,$ see \cite{stamatiou_halidias:2019}. We show that the truncated SD method is able to preserve the asymptotic stability of the underlying SDE. Motivated by a numerical example, we also propose a different SD scheme, where we actually approximate first the Lamperti transformation of the original SDE. We call this scheme Lamperti semi-discrete (LSD). It preserves positivity (in this case) of the solution, has similar asymptotic properties as the other versions of the SD method and seems promising, since there is no need for an exponential calculation. We will study the LSD method, and its properties in a forthcoming paper.

\bibliographystyle{unsrt}\baselineskip12pt
\bibliography{asSD}

\appendix

\section{Solution of linear SDEs in the narrow sense}\label{asTSD-ap:solution_linear_SDEs}

Consider the following linear in the narrow sense SDE, 
\beqq\label{asTSD:apn:linearSDE}
x_t = x_{t_0} + \int_{t_0}^t ax_s ds + \int_{t_0}^t bdW_s,
\eeqq
for $t\geq t_0,$ where $a,b$ are constants. Applying the It{\^o} formula to the transformation $U(t,x) = e^{-a(t-t_0)}x$, we obtain 

\beao
dU(t,x) & = & \left(\frac{de^{-a(t-t_0)}}{dt}x_t  + ax_te^{-a(t-t_0)}\right)dt +  b e^{-a(t-t_0)}dW_t\\
	& = & b e^{-a(t-t_0)}dW_t,
\eeao
or 
\beao
e^{-a(t-t_0)}x_t & = & x_{t_0} + b \int_{t_0}^te^{-a(s-t_0)}dW_s\\
      x_t & = & e^{a(t-t_0)}x_{t_0} + b e^{a(t-t_0)}\int_{t_0}^te^{-a(s-t_0)}dW_s.
\eeao

\section{Positivity of (\ref{asTSD-eq:exampleSDE})}\label{asTSD-ap:positivity_of_SDE}

In order to prove that $x_t>0$ a.s. we first show moment bounds of the SDE (\ref{asTSD-eq:exampleSDE}). 

\ble[Uniform moment bounds for $(x_t)$]\label{asTSD:lemma:momentbounds}
The solution process $(x_t)$ of SDE (\ref{asTSD-eq:exampleSDE}) satisfies
$$
\bfE(\sup_{0\leq t\leq T}(x_t)^p)<A,
$$
for some $A>0$ and any integer $p$ with $2\leq p \leq 19/2.$ \ele
\bpf[Proof of Lemma \ref{asTSD:lemma:momentbounds}] 
For all $|x|>R$ with $R>1,$ we have that 
\beao
J(x)&:=&\frac{xa(x) + (p-1)b^2(x)/2}{1+x^2}=\frac{x(-10x^3) + (p-1)x^{4}/2}{1+x^2}\\
&=&\frac{-21+p}{2}\frac{x^4}{1+x^2}\leq0, 
\eeao 
where the last inequality is valid for all $p$ such that $p \leq 21.$ Thus $J(x)$ is bounded for
all $x\in\bbR,$ since when $|x|\leq R.$ Application of \cite[Th. 2.4.1]{mao:2007} implies
$$
\bfE (x_t)^p\leq 2^{(p-2)/2}(1 + (x_0)^p),
$$
for any $2\leq p \leq 21,$ since $x_0\in\bbR.$ Using It\^o's formula on $(x_t)^p,$ with $p \leq 19/2$ (in order to use Doob's martingale inequality later) we have that
\beao
(x_{t})^p &=& (x_0)^p + \int_{0}^{t}\left(p(x_s)^{p-1}(-10x_s^3)  +  \frac{p(p-1)}{2}(x_s)^{p-2}x_s^{4}\right) ds\\
&& + \int_{0}^{t}p(x_s)^{p-1}x_s^2dW_s\\
&\leq& (x_0)^p + p\int_{0}^{t}(-10 + \frac{p-1}{2})(x_s)^{p+2}ds + M_t\\
&\leq& (x_0)^p + M_t,
\eeao
for any even $p$ with $2\leq p\leq21,$ or $p = 21,$ where $M_t=\int_{0}^{t}p(x_s)^{p+1}dW_s.$ Taking the supremum and then expectations in the above inequality we get
\beao
\bfE(\sup_{0\leq t\leq T}(x_{t})^p) &\leq& \bfE(x_0)^p + \bfE\sup_{0\leq t\leq T}M_t\\
&\leq&(x_0)^p + \sqrt{\bfE\sup_{0\leq t\leq T}M_t^2}\\
&\leq&(x_0)^p + \sqrt{4\bfE M_T^2},
\eeao
where in the last step we have used Doob's martingale inequality to the diffusion term $M_t.$ 
\epf

\ble[Positivity of $(x_t)$]\label{asTSD:lemma:positivity} The solution process $(x_t)$ of SDE (\ref{asTSD-eq:exampleSDE}) is positive in the sense that $x_t > 0$ a.s. 
\ele

\bpf[Proof of Lemma \ref{asTSD:lemma:positivity}] 
Set the stopping time $\theta_R=\inf\{t\in [0,T]: x_t^{-1}>R\},$ for some $R>0,$ with the convention that $\inf\emptyset=\infty.$ Application of It\^o's formula on $(x_{t\wedge\theta_R})^{-2}$ implies,
\beao
&&(x_{t\wedge \theta_R})^{-2} = (x_0)^{-2} + \int_{0}^{t\wedge \theta_R}(-2)\left((x_s)^{-3}(-10)(x_s)^3 + 3(x_s)^{-4}(x_s)^4 \right) ds\\
&& + \int_{0}^{t\wedge \theta_R}(-2)(x_s)^{-3}x_{s}^{2}dW_s\\
&\leq& (x_0)^{-2} + \int_{0}^{t\wedge \theta_R}23ds +\int_{0}^{t}(-2)x_s^{-1}dW_s\\
&\leq&(x_0)^{-2} + 23T + M_t,
\eeao
where $M_t:=\int_{0}^{t}(-2)x_s^{-1}\bbi_{(0,t\wedge\theta_R)}(s)dW_s.$
Taking expectations in the above inequality and using the fact that $\bfE M_t=0,$ we get that
$$
\bfE(x_{t\wedge \theta_R})^{-2}\leq (x_0)^{-2} + 23T< C,$$
with $C$ independent of $R.$ 
Therefore 
$$
\bfE\left(\frac{1}{x_{t\wedge\theta_R}^{2}}\right)=R^2\bfP(\theta_R\leq t) +
\bfE\left(\frac{1}{x_t^{2}}\bbi_{(t<\theta_R)}\right)<C,
$$
implying that
$$
\bfP(x_t\leq0)=\bfP\left(\bigcap_{R=1}^\infty \Big\{x_t<
\frac{1}{R} \Big\}\right)=\lim_{\Rto}\bfP\left( \Big\{x_t<
\frac{1}{R} \Big\}\right)\leq  \lim_{\Rto}\bfP(\theta_R\leq t)=0.
$$
We conclude that $x_t > 0$ a.s. 
\epf

\section{Lamperti Tranformation of (\ref{asTSD-eq:exampleSDE})}\label{asTSD-ap:Lamperti_tranformation}

Applying the It\^o formula to the transformation $z(x) = -1/x,$ we obtain 

\beao
dz_t & = & \left( (x_t)^{-2}(-10)(x_t)^3  + \frac{1}{2}(-2)(x_t)^{-3}(x_t)^4\right)dt +  (x_t)^{-2}(x_t)^{2}dW_t\\
& = &  -11x_t + dW_t\\
& = &  11(z_t)^{-1}dt + dW_t\\
\eeao
or for $t\geq t_0$
\beao
z_t & = & z_{t_0} + 11 \int_{t_0}^t (z_s)^{-1} + \int_{t_0}^t dW_s\\
 & = & z_{t_0} + 11 \int_{t_0}^t (z_s)^{-1} + W_t - W_{t_0}.
\eeao

\section{Stochastic Integral Approximation}\label{asTSD-ap:integralapprox}

We want to estimate the stochastic integral appearing in the proposed truncated semi-discrete method (\ref{asTSD-eq:exampleSDtrunc}) for the approximation of SDE (\ref{asTSD-eq:exampleSDE}). In a similar way we calculate the integral appearing in the exponential truncated semi-discrete scheme (\ref{asTSD-eq:exampleSD}). 

In the numerical simulations we used the following relation
$$\int_{t_n}^{t_{n+1}}e^{10\pi^2_\D(y_{n}^\D)s}dW_s\approx e^{10\pi^2_\D(y_{n}^\D)t_n}\D W_n.$$
We show the following estimation
\beqq\label{asTSD:apn:maximalmartingale}
\bfP\left(\left|\int_{t_n}^{t_{n+1}}e^{10\pi^2_\D(y_{n}^\D)s}dW_s-e^{10\pi^2_\D(y_{n}^\D)t_n}\D W_n\right|\geq\D^r\right)\leq2e^{20\pi^2_\D(y_{n}^\D)t_{n+1}}\D^{1-2r},
\eeqq
suggesting that the probability of the absolute difference of these two random variables being of order $\D^r,$ with $0<r<1/2,$ approaches unity as $\D$ goes to zero.
First, we write the difference of the two local martingales as 
$$
\int_{t_n}^{t_{n+1}}e^{10\pi^2_\D(y_{n}^\D)s}dW_s - e^{10\pi^2_\D(y_{n}^\D)t_n}\D W_n = \int_{t_n}^{t_{n+1}}\left(e^{10\pi^2_\D(y_{n}^\D)s} - e^{10\pi^2_\D(y_{n}^\D)t_n}\right)dW_s
$$
and then use the martingale inequality to get for any $\ep>0$ that
\beao
&&\bfP\left(\left|\int_{t_n}^{t_{n+1}}\left(e^{10\pi^2_\D(y_{n}^\D)s} - e^{10\pi^2_\D(y_{n}^\D)t_n}\right)dW_s\right|\geq \ep\right) \\
&\leq&\ep^{-2}\bfE\left(\left(\int_{t_n}^{t_{n+1}}\left(e^{10\pi^2_\D(y_{n}^\D)s} - e^{10\pi^2_\D(y_{n}^\D)t_n}\right)dW_s\right)^2 \bigg|\bbf_{t_n}\right)\\
&\leq&\ep^{-2}\int_{t_n}^{t_{n+1}}\bfE\left(\left(e^{10\pi^2_\D(y_{n}^\D)s} - e^{10\pi^2_\D(y_{n}^\D)t_n} \right)^2\bigg|\bbf_{t_n}\right)ds\\
&\leq&\ep^{-2}\left[\frac{1}{20\pi^2_\D(y_{n}^\D)}\left(e^{20\pi^2_\D(y_{n}^\D)t_{n+1}} - e^{20\pi^2_\D(y_{n}^\D)t_n}\right) + e^{20\pi^2_\D(y_{n}^\D)t_{n}}\D\right]\\
&& - 2\ep^{-2}e^{10\pi^2_\D(y_{n}^\D)t_{n}}\frac{1}{10\pi^2_\D(y_{n}^\D)}\left(e^{10\pi^2_\D(y_{n}^\D)t_{n+1}} - e^{10\pi^2_\D(y_{n}^\D)t_n}\right)\\
&\leq&\ep^{-2}e^{20\pi^2_\D(y_{n}^\D)t_{n+1}}\left(\frac{1}{20\pi^2_\D(y_{n}^\D)}\left(1 - e^{-20\pi^2_\D(y_{n}^\D)\D}\right) + e^{-20\pi^2_\D(y_{n}^\D)\D}\D\right)\\
&\leq&2\ep^{-2}e^{20\pi^2_\D(y_{n}^\D)t_{n+1}}\D,
\eeao
where in the last step we used the inequality $1 - e^{-x}\leq x,$ for any $x>0.$ We apply the above inequality for $\ep = \D^r,$ with $0<r<1/2$ to get (\ref{asTSD:apn:maximalmartingale}).

\end{document}